\documentclass[10pt]{article}
\usepackage{amssymb}
\font\elevenss=cmss11

\font\eightss=cmss8

\font\sixss=cmss8 at 6pt

\newfam\ssfam
\textfont\ssfam=\elevenss \scriptfont\ssfam=\eightss
 \scriptscriptfont\ssfam=\sixss
\newtheorem {thm}{Theorem}[section]

\newtheorem {pr}[thm]{Proposition}

\def\Cox{\hfill \square}

\def\grad{\bigtriangledown}

\def\ee{\epsilon}

\def\cc{{\vec c}}
\def\hh{{\vec h}}

\def\W{{\vec W}}

\def\E{{\mathbb{E}}}

\def\P{{\mathbb{P}}}

\def\F{{\cal{F}}}

\def\|{\, | \, }
\def\one{{\bf 1}}

\def\Qt{\tilde{Q}}
\def\Pt{\tilde{\P}}
\def\nbd{{\cal N}}

\begin{document}

\begin{titlepage}
\begin{center}
{\large \bf Time to absorption in discounted reinforcement models} \\
\end{center}
\vspace{5ex}
\begin{flushright}
Robin Pemantle \footnote{Research supported in part by National
Science Foundation grant \# DMS 0103635}$^,$\footnote{Department
of Mathematics, The Ohio State University, 231 W. 18 Ave.,
Columbus, OH 43210}
  ~\\
Brian Skyrms \footnote{ School of Social Sciences, University of
California at Irvine, Irvine, CA 92607} \\

\end{flushright}

\vfill

{\bf ABSTRACT:} \hfill \break
Reinforcement schemes are a class of non-Markovian stochastic processes.
Their non-Markovian nature allows them to model some kind of memory 
of the past.  One subclass of such models are those in which the past is
exponentially discounted or forgotten.   Often, models in this subclass
have the property of becoming trapped with probability~1 in some 
degenerate state.  While previous work has concentrated on such limit
results, we concentrate here on a contrary effect, namely that the
time to become trapped may increase exponentially in $1/x$ as the
discount rate, $1-x$, approaches~1.  As a result, the time to become
trapped may easily exceed the lifetime of the simulation or of the
physical data being modeled.  In such a case, the quasi-stationary
behavior is more germane.  We apply our results to a model of social network
formation based on ternary (three-person) interactions with uniform
positive reinforcement.
\vfill


\end{titlepage}

\setcounter{equation}{0}
\section{Introduction} \label{ss:intro}

Stochastic models incorporating mechanisms by which likelihoods of
outcomes increase according to their accumulated frequencies have
been around since the introduction of P\'olya's Urn~\cite{EP23}.
The mathematical framework for many of these models appears in
the literature on stochastic approximation, beginning with~\cite{RM51},
in the urn model literature~\cite{Fr65}, in the literature on
reinforced random walks~\cite{Pem90,Dav90,Lim01}, and in the
literature on the relation between stochastic systems and their
deterministic mean-field dynamical system approximations~\cite{BH95,Ben96}.  

These processes, known in the mathematical community as 
{\em reinforcement processes}, have long been used by psychologists 
as models for learning~\cite{BM55,IT69,Nor72,Lak81}.  
Increasingly, reinforcement models have been adopted by other social 
scientists as interactive models in which collective learning takes 
place in the form of network formation or adaptation of strategies:
sociologists studying the ``small world'' network phenomenon~\cite{WS98,BW00}, 
formation of dyads of reciprocal approval~\cite{FM96}; 
economists studying evolutionary game theory~\cite{MS82},
strategic learning~\cite{RE95}~\cite{FK93} or its interaction with network 
structure~\cite{Ell93,AI97}.  These models are designed to 
explore mechanisms by which agents with limited information, rationality 
or sophistication may nevertheless achieve advantageous social structures 
via the application of simple rules for behavior change.

Due, perhaps, to a dearth of types of simple local rules (or perhaps
to a lack of imagination on the part of modelers), most reinforcement
models fall into one of two classes.  The first class contains models
for which the past is weighted uniformly.  This class includes the
urn models, stochastic approximations and reinforced random walks 
mentioned above, as well as the economic game theory models of Erev 
and Roth~\cite{RE95}.  Uniform weighting means that the step from time
$n-1$ to $n$ represents a fraction of only $1/n$ of the total learning
up to time $n$, so one obtains a time-inhomogenous process in which
the hidden variables change by amounts on the order of $1/n$ at time $n$.
The second class consists of models in which the past is 
exponentially discounted or forgotten.  This class includes
the learning models of the 1960's and 70's~\cite{IT69}, as well as many 
contemporary models of repeated economic games, e.g.~\cite{BL02}.  
In these models, the weight of the present is asymptotically equal
to the discount parameter, $x$, defined as the $x$ for which 
an action $t$ units of time will be weighted by $(1-x)^t$.  
More precisely, the fraction of total learning between time $n-1$ and $n$ 
of the total learning to time $n$ will be roughly the maximum of $1/n$
and $x$.  

Our chief concern in this paper is to study how the discounted
process approaches the non-discounted process as $x \to 0$.  The
long-run behavior in these two cases is qualitatively different.
Limit theorems for non-discounted processes have been obtained chiefly
in the framework of stochastic perturbations to dynamical systems.
Typically, the stochastic system converges to limit points or limit
cycles of the dynamical system that corresponds to the mean motion
of the stochastic system~\cite{BH95}.  The random limit is supported on
weakly stable equilibria~\cite{Pem90}, though the system may remain
near unstable equilibria for long periods of time (see~\cite{BP03}
for a discussion of this phenomenon in continuous time, and~\cite{PS03b}
for a case study in discrete time).  

In the discounted processes we study here, there are trapping states, 
into which the chain must eventually fall (there is another kind of
discounted process we are not concerned with here, which converges to 
an ergodic Markov chain, see, e.g.,~\cite{IT69}).  The reciprocity model
of~\cite{BL02}, for example, is of this type.  From the point of view
of studying the transition as discounting goes to zero, the most 
interesting case is when the trapping states are disjoint from 
the stable equilibria in the non-discounted process.  Trapping states for 
the discounted process must always be equilibria for the non-discounted
process, but when all the trapping states are unstable equilibria 
(equivalently, all the stable equilibria of the non-discounted process
are non-trapping in the discounted process), then the conflict between
the discounted and non-discounted behavior is maximized.

The transition is easy to describe informally.  As the discount rate
approaches zero, the discounted process behaves for a longer and longer
time like the non-discounted process, and then abruptly falls into a trap.
Of course when the discounting parameter is $x$ it takes time 
at least on the order of $1/x$ for the system to notice there is
discounting going on.  But in fact, due to the learning that has
gone on during this phase, it will take time of order $\exp (c x^{-1})$
before the system discovers a trap and falls in.  It is in fact not 
hard to guess this via back-of-the-napkin computations.
One of our main motivations for pursuing this rigorously was to explain 
why simulation data contradicted the easily proved limit theorem: it
was because the time scale of the simulation (let alone of any real 
phenomenon modeled by the simulation) was never anywhere near the
time needed to find a trap.  

Our purpose in the present paper is to prove various versions of this.
In the next section, we present the ternary interaction model which
was our original motivation for this study.  Section~\ref{ss:one-d}
then introduces a simple process that is a building block for the
ternary interaction model.  For that process, results about trapping times
can be proved with the correct constant.  The last section then
proves a $\exp (c x^{-1})$ waiting time result for a general class of
models, but without the correct value of $c$.  With a little
linear algebra, this is shown to apply to the ternary interaction model
of Section~\ref{ss:three's company}.

\setcounter{equation}{0}
\section{Three's Company: a ternary interaction model} 
\label{ss:three's company}

The following process is described in~\cite{PS03a}, where it is called
{\em Three's Company}, and is put forth as a model
for formation of ternary collaborations in a three-player
version of Rousseau's stag hunting game.
Fix a positive integer $N \geq 4$, representing the size of the 
population.  For $t \geq 0$ and $1 \leq i , j \leq N$, define 
random variables $W(i,j,t)$ and $U(i,t)$ inductively
on a common probability space $(\Omega , \F , \P)$ as follows.
The $W$ variables are positive numbers, and the $U$ variables are
subsets of the population of cardinality~3.  One may think of the
$U$ variables as random triangles in the complete graph with a 
vertex representing each agent.  The initialization is
$W(i,j,0) = 1$ for all $i \neq j$, while $W(i,i,0) = 0)$.  The
inductive step, for $t \geq 0$, defines probabilities (formally, 
conditional probabilities given the past) for the variables $U(i,t)$ 
in terms of the variables $W(r,s,t)$, $r , s \leq N$, and then defines 
$W(i,j,t+1)$ in terms of $W(i,j,t)$ and the variables $U(r,t)$, $r \leq N$.
The equations are:
\begin{eqnarray}
\P (U(i,t) = S \| \F_t) & = & \frac{\one_{i \in S} \prod_{r,s \in S, r < s}
   W(r,s,t)}{\sum_{S' : i \in S'} \prod_{r,s \in S' , r < s}
   W(r,s,t)} \, ; \label{eq:probs} \\[1ex]
W(i,j,t+1) & = & (1-x) W(i,j,t) + \sum_{r=1}^N \one_{i,j \in U(r,t)} \, .
   \label{eq:updates}
\end{eqnarray}
Here $(1-x)$ is the factor per unit time by which the past is discounted,
and the $\sigma$-field conditioned on is the process up to time $t$,
$$\F_t := \sigma \left \{ W(i,j,u) : u \leq t \right \} \, .$$
We may think of the normalized matrix 
$${\bf W}_t := \frac{1}{\sum_{i,j} W(i,j,t)} \, W(\cdot , \cdot , t)$$ 
as the state vector, which is then an asymptotically time-homogeneous 
Markov chain, with an evolution rule of the well known form
\begin{equation} \label{eq:state}
\E \left ( {\bf W}_{t+1} - {\bf W}_t \| \F_t \right ) = 
   g(t) \left [ \mu ({\bf W}_t) + \xi_t \right ] 
\end{equation}
where in this case, $g(t) = 1/x + O(1/t)$, the drift vector field $\mu$ 
maps the simplex of normalized matrices into its tangent space and may 
be explicitly computed, and $\xi_t$ are martingale increments of order~1.

These equations model a social interaction in which each agent $i$ at each 
time $t$ invites two others to frolic\footnote{Engage in some rewarding 
interaction such as anti-competitive price fixing}.  For each agent $i$,
the trio chosen by $i$ is chosen from all possible trios containing $i$, 
according to the products of the weights $W(i, \cdot , t)$.  
Thus the probability of agent $i$ forming the trio $\{ i , j , k \}$ is 
proportional to $W(i,j,t) W(i,k,t)$.  After the frolicking, fond memories 
ensue: each of the three pair weights $W(i,j,t), W(i,k,t)$ {\em and} 
$W(j,k,t)$ is increased by~1, and a portion $x$ of the past weights 
is forgotten.  For ease of bookkeeping, the weights of unordered
pairs are defined as symmetric weights of ordered pairs, so the 
weights $W(j,i,t) , W(k,i,t)$ and $W(k,j,t)$ are increased as well.
We write $W(e,t)$ for $W(i,j,t)$ when $e$ is the edge (unordered set)
$\{ i , j \}$.  

It is shown in~\cite{PS03a} that the network always breaks into small
cliques, with interactions occurring only among the cliques:
\begin{thm} \label{th:limit}
In Three's Company, with population size $n \geq 4$ and any
discount rate $x \in (0,1)$, with probability~1 the population may 
be partitioned into subsets of sizes~3, 4 and~5, such that each member 
of each subset chooses each other with positive limiting frequency, 
and chooses members outside the subset only finitely often.  Every partition 
into sets of sizes~3, 4 and~5 has positive probability of occurring.  $\Cox$
\end{thm}

Simulation data is also given there.  For $N=6$, if $x = .4$ (a rather
steep discount rate), the network always breaks into two cliques of size~3,
as predicted by the theorem.  When $N = 6$ and $x = .2$, which is still
a greater discount rate than one finds in most economic models, one
finds, with runs of several thousand, that no such structure emerges.
Instead, all six members of the population remain well connected.  
This is because of the exponential time scale of the transition from
stable equilibria (well-connectedness is a stable equilibrium of the
non-discounted model) to trapping states (two cliques of size three is
the unique trapping state of the $N=6$ discounted model).  Specifically,
in the last section of this paper we will prove:
\begin{thm} \label{th:transience time}
In the game Three's Company, for each $N \geq 6$ there is a $\delta > 0$ 
and numbers $c_N > 0$ such that in Three's Company with $N$ players 
and discount rate $1-x$, the probability is at least $\delta$ that each 
player will play with each other player beyond time $\exp (c_N x^{-1})$.
\end{thm}

\setcounter{equation}{0}
\section{Trapping in one-dimensional discounted reinforcement} 
\label{ss:one-d}

In this section we analyze a one-dimensional process in which 
sharp quantitative results may be obtained on the exponential
rate at which the time until trapping increases with $1/x$,
where $1-x$ is the discount factor.  
This is in keeping with our philosophy of providing sharp results 
on a collection of simplified models that constitute building
blocks for more complicated models.  In the last section we will
apply the principles gleaned from this to get bounds on the 
exponential rate of increase of trapping time in the Three's Company model.

Let us consider a system whose state vector
varies in the interval $[0,1]$ with evolution dynamics that are
symmetric around an attractor at $1/2$, and whose transitions from
state $w$ have a profile that depends on $w$ and is scaled by $x$.  
In analogy with models such as Three's Company, we assume that 
the unscaled transitions have variance bounded from below as $w$ 
varies over compact sub-intervals of $(0,1)$.  Thus the rules of 
evolution of the state vector $W$ may be given in terms of 
probability distributions $Q_w$, parametrized by $w \in [0,1]$, 
with bounded support, satisfying $Q_w (s) = Q_{1-w} (-s)$, and obeying:
\begin{equation} \label{eq:motion}
\P (W(n+1) - W(n) \in x \cdot S \| \F_n) = Q_{W(n)} (S)
\end{equation} 
on the event that $W(n)$ is in a compact subinterval $I_x$, with
$I_x \uparrow (0,1)$ as $x \to 0$.  We assume that the mean of
$Q_w$ is positive on $(0,1/2)$ and negative on $(1/2 , 1)$, but that
$Q_w$ has both positive and negative elements in its support and
varies smoothly with $w$.

As an example, one may consider a class of two-color urn 
models generalizing Friedman's Urn~\cite{Fr49,Fr65} in the discounted 
setting.  An urn begins with $R(0)$ red balls and $B(0)$ black balls.  
At the $n^{th}$ time step, a random number $U(n)$ of red balls and
$V(n)$ black balls are added to the urn.  Conditional on the past,
$\P (U(n) = k) = u_{W(n-1)} (k)$ and $\P (V(n) = k) = u_{1-W(n-1)} (k)$, 
where $W (n) := R(n) / (R(n) + B(n))$ is the state parameter, in 
this case the proportion of red balls, and $u_w$ are probability 
distributions on the nonnegative integers, continuously varying 
in the parameter $w \in [0,1]$, satisfying 
$$\frac{\sum_k k u_w (k)}{\sum_k k (u_w (k) + u_{1-w} (k))} > w
   \mbox{ for } 0 < w < \frac{1}{2} \, .$$
At the end of each step, all balls are reduced in weight by a factor 
of $1-x$.   For greater specificity, one may keep in mind an example 
where two balls are sampled: if they are of the same color then one 
ball of that color is added; if they are of different colors then 
one ball of each color is added.  

In the non-discounted system, where the step size scales as $1/n$ at
time $n$ instead of holding constant at $x$, the system is well 
approximated by a diffusion with incremental variance of order $n^{-2}$ 
and drift $n^{-1} \mu_w$, with $\mu_w$ being the mean, $\overline{Q_w}$,
of $Q_w$.  Thus~(\ref{eq:state}) holds with $g(t) = t^{-1}$ and $\mu (w) = 
\overline{Q_w}$.  The system must converge to the unique attracting 
equilibrium at $1/2$~\cite{Pem90}.  In the discounted case, although the 
state must converge to~0 or~1, the logarithm of the expected time
to come near~0 or~1 may be computed in terms of the following data.

Pick any $w \in (0,1/2)$.  The quantity 
$$Z_w (\lambda) := \int \exp (- \lambda y) \, dQ_w (y)$$
is equal to~1 at $\lambda = 0$.  The derivative $(d / d\lambda) 
Z_w (\lambda) |_{\lambda = 0}$ is given by 
$\int (-y) \, dQ_w (y)$, which is negative by the assumption
that $Q_w$ has positive mean.  On the other hand, since $Q_w$ gives 
positive probability to negative values, we see that as 
$\lambda \to \infty$, $Z_w \to \infty$, and by convexity of $Z_w (\cdot)$ 
it follows that there is a unique $\lambda_w > 0$ for which 
$Z_w (\lambda_w) = 1$.  Define
$$\Lambda (w) := \int_w^{1/2} \lambda_u \, du$$
and let $C := \Lambda (0)$.

\begin{thm} \label{th:rate}
Let $I_x \uparrow (0,1)$ as $1-x \uparrow 1$, slowly enough so that
transitions outside of $[0,1]$ are never possible.  Let $T_x$ be the
expectation of the first time $n$ that  $W(n) \notin I_x$.  Then
as $1-x \uparrow 1$,
$$x \log \E T_x \to C \, .$$
\end{thm}

\noindent{\bf Remark:} This is essentially a large deviation problem,
so the rate $C$ is not determined by the mean and variance of $Q_w$
but rather by the exponential moments of $Q_w$.  In particular, there
are many processes which satisfy~(\ref{eq:state}) with the same
$g$ and $\mu$, but their large deviation rates depend on the fine structure
of the increment distribution through the exponential moments, as 
captured by $Z_w$ and $\lambda_w$.  The solution of this rate problem
is standard; a similar analysis may be found, for example
in~\cite[Section~5.8.2]{DZ93}.  

\noindent{\sc Proof:}
For one inequality, we fix any $\delta > 0$.  Define the quantity 
$$M_{(\delta)} (t) := \exp ((1 - \delta) x^{-1} \Lambda (W(t))  \, .$$
Since $Q_w$ varies smoothly with $w$ and has bounded support, 
we see that as $x \to 0$, 
$$\Lambda (W(t+1)) - \Lambda (W(t)) = (W(t+1) - W(t)) (- \lambda_{W(t)} 
   + O(x)) \, .$$
Therefore, since the conditional distribution of $x^{-1} (W(t+1) - W(t))$
given $\F_t$ is given by $Q_{W(t)}$, we see that
\begin{eqnarray}
\E (M_{(\delta)} (t+1) \| \F_t) & = & M_{(\delta)} (t) \, \E \left [ 
   \exp ((1 - \delta) x^{-1} (W(t+1) - W(t)) (- \lambda_{W(t)} + O(x))) 
   \right ] \nonumber \\
& \to & M_{(\delta )} (t) Z_{W(t)} ((1 - \delta) \lambda_{W(t)}) 
   \label{eq:super}
\end{eqnarray}
uniformly in $W(t)$ as $x \to 0$.  We know that $Z_w (\cdot ) < 1$ on 
$(0 , \lambda_w)$, hence we may pick $x = x(\delta)$ small enough so that 
$$M_{(\delta)} (t)^{-1} \E (M_{(\delta)} (t+1) \| \F_t) < 1$$
or in other words, so that $M_{(\delta)}$ is a supermartingale.

Let $I_x = [a_x , 1 - a_x]$.
Starting with $W(t_0) \in (1/2 - \delta , 1/2)$, and stopping at
the time $\tau$ when $W(\cdot)$ exits $[a_x,1/2]$, we have
for some constant $c(\delta)$ going to zero with $\delta$,
\begin{eqnarray*}
\exp (x^{-1} c(\delta)) & \geq & M_{(\delta)} (t_0) \\[1ex]
& \geq & \E ( M_{(\delta)} (\tau) \| \F_{t_0} ) \\[1ex]
& \geq & \P ( M_{(\delta)} (\tau) < a_x ) \exp (x^{-1} \Lambda (a_x) )
\end{eqnarray*}
which implies that 
$$\log \P ( M_{(\delta)} (\tau) < a_x ) < - (1 - \delta) x^{-1} 
   (\Lambda (0) + o(1))$$
as $x \to 0$.  A completely analogous argument shows that
the process started in $(1/2 , 1/2 + \delta)$ exits $[1/2 , 1- a_x]$
with at most this probability as well.  The trajectory of $W(\cdot)$
may be decomposed into segments that begin in $[1/2 - \delta , 1/2 + \delta]$
and end when $W(t) - 1/2$ changes sign or $W(t) \notin I_x$.  We have shown
that the expected number of trajectories is at least
$\exp (x^{-1} (1 - \delta) (\Lambda (0) + o(1)))$ as $x \to 0$, 
which implies that the number of time steps until exiting $I_x$ is 
at least this great, once $W(t) \in [1/2 - \delta , 1/2 + \delta]$.  
Letting $c' (\delta)$ denote the probability of entering this interval,
we see that 
$$\E T_x \geq c' (\delta) \exp (x^{-1} (1 - \delta) (\Lambda (0) + o(1))$$
as $x \to 0$, and finally, sending $\delta$ to zero proves that
$$\liminf x \log \E T_x \geq C \, .$$

For the other direction, define a tilted measure on the space
of trajectories $\{ W(t) : t = 0 , 1 , 2 , \ldots  \}$ as follows.
The equation~(\ref{eq:motion}) is replaced by
$$\Pt (W(n+1) - W(n) \in x \cdot S \| \F_n) = \Qt_{W(n)} (S)$$
where $\delta > 0$ is fixed and the Radon-Nikodym derivative is given by
$$\frac{d\Qt_w}{dQ_w} (y) 
   = \frac{\exp \left [ (1+\delta) (\Lambda (w+y) - \Lambda (w)) \right ]}
     {\int \exp \left [ (1+\delta) (\Lambda (w+y) - \Lambda(w)) \right ] 
     \, dQ_w (y) } \, .$$

The measure $\Pt$ is designed to have two properties.  First, the process
$\{ W(t) \}$ is a supermartingale on $[a_x , 1/2]$ with respect to 
$\Pt$ for sufficiently small $x$.  To see this, note that this is 
equivalent to $\Qt_w$ having negative mean, which is equivalent to
$$\int y e^{(1 + \delta) (\Lambda (w+y) - \Lambda (w))} \, dQ_w (y) \leq 0$$
for all $w \in [a_x , 1/2]$.  The quantity $\Lambda (w + y) - \Lambda (w)$ 
is equal to $y (- \lambda_w + o(1))$ as $x \to 0$, so it suffices to show
that 
$$\int y e^{-(1 + \delta) y \lambda_w} \, dQ_w (y) < 0 \, .$$
But this follows from the fact that $Z_w$ is convex and increases 
through~1 at $\lambda_w$: the derivative at $(1 + \delta) \lambda_w$
must therefore be positive, and the derivative may be identified as
$$- \int y e^{- (1 + \delta) \lambda_w y} \, dQ_w (y) \, ,$$
proving the supermartingale property.

The second property is that if $\tau$ is the exit time of 
$[a_x , 1/2]$, then on the $\sigma$-field $\F_\tau$, 
$d \Pt / d\P$ is at most $\exp ((1+\delta)x^{-1} \Lambda (0))$.  
Indeed, by its definition, 
$$\frac{d\Pt}{d\P} (W(t_0) , \ldots , W(\tau)) = \prod_{t=t_0}^{\tau - 1}
   \frac{\exp \left [ x^{-1} (1+\delta) (\Lambda (W(t+1) - \Lambda (W(t))) 
   \right ]}{\int \exp \left [ (1+\delta) (\Lambda (W(t)+y) - \Lambda(W(t))) 
   \right ] \, dQ_{W(t)} (y) } \, .$$
The denominator of each factor is at least~1 by the fact that
$$M_{(-\delta)} (t) := \exp ((1 + \delta) x^{-1} \Lambda (W(t)))$$
is a submartingale when $x(\delta)$ is small enough, which is proved
by a computation exactly analogous to~(\ref{eq:super}).  The product of
the numerators is simply $\exp ((1+\delta) x^{-1} \Lambda (W(\tau)) - 
\Lambda (W(t_0)))$, which is at most $\exp ((1+\delta) x^{-1} \Lambda (0))$, 
proving the second property.

Running $\Pt$ on $[a_x , 1/2]$ and its reflection on $[1/2 , 1-a_x]$,
the process $1/2 - |1/2 - W(t)|$ is a supermartingale with incremental 
variance of order $x^{-2}$.  The median time for it to reach a value
less than $a_x$ is therefore at most $O(x^{-2})$.  Comparing $\Pt$ and
$\P$, we find that there is a $c$ such that from any starting data, 
the probability of exiting $I_x$ by time $c x^{-2}$ is at least
$(1/2) (d\P / d\Pt) \geq (1/2) \exp (- (1+\delta)x^{-1} \Lambda (0))$.  
Breaking into time intervals of size $c x^{-2}$, it then follows 
that the mean time to exit $I_x$ is at most 
$2 c x^{-2} \exp (C (1+\delta) x^{-1})$.  As this holds for any
$\delta > 0$ (and constants depending on $\delta$), this proves that
$$\limsup x \log \E T_x \leq C$$
and finishes the proof of the theorem.   $\Cox$

\setcounter{equation}{0}
\section{Proof of Theorem~\protect{\ref{th:transience time}}} 
\label{ss:proof}

In analogy with the one-dimensional toy model, we expect to find an
exponential wait to trapping in Three's Company if the non-discounted 
system has an attractor outside of the limit set of absorbing states of 
the discounted system.  Unfortunately, at this point we cannot see 
any way to compute the large deviation rate in multi-dimensional problems.
The standard multi-dimensional analogue to Theorem~\ref{th:rate} 
is expressed as a variational result involving minimizing a functional 
over all paths.  We settle for proving the existence of a nonzero 
exponential rate in $x^{-1}$.  The following result will imply
Theorem~\ref{th:transience time}.

\begin{pr} \label{pr:multi}
Let the vector-valued Markov chain $\W (t)$ satisfy
\begin{equation} \label{eq:dynamics}
\P (\W(t+1) - \W(t) \in x S \| \W(t)) = Q_{\W(t)} (S)
\end{equation}
with $Q_w$ having bounded support and varying smoothly as $w$ 
varies over some closed neighborhood $\nbd$ of a point $\cc$.  
Suppose there is a strong Lypunov function $V$, meaning that $V$ is 
smooth and bounded with $V < 0$ on $\nbd^c$, $V(\cc) > 0$ and
$$\int V(w + y) \, dQ_w (y) > V(w)$$ 
for all $w \in \nbd$.  Then there is a constant $\gamma$ such that
for all $\hh$ in some smaller neighborhood of $\cc$, 
$$\E_\hh T_x > \gamma \exp (\gamma x^{-1})$$
for sufficiently small $x$, with $T_x$ being the time to exit $\nbd$.
\end{pr}

\noindent{\sc Proof:} Given a non-negative parameter, $\lambda$, 
define $M(t) = M_\lambda (t) = \exp (- \lambda V(\W(t)))$.  Arguing
as in the first half of the proof of Theorem~\ref{th:rate}, we see
from the bounded support hypothesis that for fixed $w \in \nbd$, 
\begin{equation} \label{eq:lam}
M_\lambda (t)^{-1} \E (M_\lambda (t+1) - M_\lambda (t) \| \F_t)
\end{equation}
vanishes at $\lambda = 0$ and has negative derivative.  By compactness
of $\nbd$ and smoothness of $Q_w$, we may choose a $\lambda > 0$ so 
that~(\ref{eq:lam}) is negative for all $w \in \nbd$, implying that
$M (t)$ is a supermartingale up to the exit time of $\nbd$.  
Let $\nbd'$ be the neighborhood $\{ \hh \in \nbd : V(h) > V(\cc) / 2 \}$.
Using $T_G$ to denote the first time $\tau \geq 0$ that $W(\tau) \in G$,
we then have, for $V(\hh) > V(\cc) / 4$,  
\begin{eqnarray*}
e^{-\lambda V(\cc) / 4} & \geq & \E_\hh (\W (0)) \\[1ex] 
& \geq & \E_\hh T_{\nbd^c \cup \nbd'} \\[1ex]
& \geq & \P_\hh (T_{\nbd^c} < T_{\nbd'})  \, . 
\end{eqnarray*}
Breaking the time $T_{\nbd^c}$ into sojourns away from $\nbd'$ 
then proves the theorem with $\gamma = \lambda V(\cc) / 4$.   $\Cox$

In Three's Company, if we start with the sum of the weights
equal to $3 N x^{-1}$ (that is, in stationarity), then the dynamics
are described exactly by~(\ref{eq:dynamics}).  We need only check the
existence of a strong Lyapunov function.  This will follow if we
can identify a hyperbolic attractor for the vector field $F(\cdot)$, where
$F(w)$ is the mean of $Q_w$.  Indeed, if $F$ vanishes at a point $\cc$ and
$dF (\cc)$ has eigenvalues with negative real parts, then there is a
quadratic function $V$ near $\cc$ satisfying $\grad V \cdot F > 0$
which we may take as the Lyapunov function.  All that remains is to 
identify the hyperbolic attractor for the mean motion field.

The mean motion is given by a vector field $F$ on the state space. 
The state space is the set of non-negative real functions $X$ on the
edges summing to $3 N x^{-1}$, which we think of as embedded in the
cone of non-negative functions, since $F$ extends naturally via
$F (\lambda X) = F(X)$.  The computations are a little more convenient
when we normalize the sum of weights to be ${N \choose 2}$.
It is also convenient to let $n = N-1$ be one less than the number of 
agents.  The attractor on which we focus is the symmetric point $\cc$
defined by $c(e) = 1$ for all $e$.  It is immediate to verify 
that $F(\cc) = 0$.  In order to verify that $\cc$ is an attractor for $F$,
we need to compute the differential of $F$ at $\cc$.  Accordingly,
let $\one_e$ denote the function that is~1 on $e$ and~0 elsewhere.
The derivative of $F$ in the $\one_e$ direction is computed as follows.

Let the edge weights at time $t$ be given by $1 + \ee \one_e$.  
The expected number of $i$ for which $f \in U(i,t)$ is $\frac{6}{n} + O(\ee)$
for all $f$.  By symmetry, the $O(\ee)$ term depends only on whether
$f$ shares two, one or zero endpoints with $e$.  For example, in
the case $f = e$, we compute the expected number of times $e$ is
reinforced as follows.  Let $e = \{ v , w \}$.  Then 
$$\P (e \in U(v,t)) = (n-1) \frac{1+\ee}{{n \choose 2} + (n-1) \ee} \, .$$
The probability of $e \in U(w,t)$ is the same.  For $z \neq v,w$, the
probability of $e \in U(z,t)$ is exactly ${n \choose 2}^{-1}$.  Summing
yields an expected increment in $W(e,t)$ of 
$$2 (n-1) \frac{1+\ee}{{n \choose 2} + (n-1) \ee} + \frac{n-1}{{n \choose 2}}
   = \frac{6}{n} + \frac{4(n-2)}{n^2} \ee + O(\ee)^2 \, .$$
We write this as
$$ \frac{6}{n} \left ( 1 + \frac{2 (n-2)}{3 n} \ee + O(\ee^2) \right ) \, .$$
Computing, the other two expectations in this manner, we find
that the expectation for $f$ at distance $j$ from $e$ is
$\frac{6}{n} (1 + B_j \ee + O(\ee^2))$, where 
\begin{eqnarray*}
B_0 & = & \frac{2 (n-2)}{3 n}  \\
B_1 & = & 0 \\
B_2 & = & - \frac{4}{3 n (n-1)} \, .
\end{eqnarray*}
From this it follows that for any $\hh$, 
$$\frac{n}{6} \E (W(t+1) \| W(t) = \cc + \ee \hh ) = \cc + M \hh 
   + O(|\hh|^2)$$
where $M$ is a generalized circulant matrix (symmetric under the 
action of edge permutation on pairs of edges) with entries $B_0$ 
on the diagonal, $B_2$ for disjoint edges, and zero otherwise.
Since $F$ is the vector field pointing toward $\E (W(t+1) - W(t) \| W(t))$,
the differential of $F$ is, up to the constant multiple $\frac{6}{n}$,
equal to $M - I$.

The eigenvalues of a matrix such as $M$ are particularly easy to 
evaluate, using the rubric of {\em association schemes}
(see~\cite[Section~2.2]{Ter96}, which is taken from~\cite{BI84,BCN89}).
All such matrices are elements of the Bose-Mesener algebra ${\cal M}$
which, in the case of the incidence graph for edges of the complete 
graph, is commutative semi-simple of dimension~3.  This implies that
$M$ has at most three distinct eigenvalues.  These may be found by 
computing the action of $M$ on the three shared eigenspaces common
to all elements of ${\cal M}$.  

The null eigenspace has dimension~1: $M \hh = 0$ if and only if 
$\hh = \lambda \cc$.  The other two eigenvalues may be gotten by 
choosing an edge $e$ and setting the eigenvectors equal to 
$a_2 H_2 + a_1 H_1 + a_0 H_0$, where $H_2 = \one_e$, $H_1$ is 
the sum of $\one_f$ over edges $f$ sharing one vertex with $e$, 
and $H_0$ is the sum of $\one_f$ over edges $f$ disjoint from $e$.  
The action of $M$ on such a sum produces another such sum, and is linear, 
having matrix
$$\frac{4}{3 n (n-1)} \left ( \begin{array}{ccc}
   {n-1 \choose 2} & 0 & -1 \\
   0 & {n-2 \choose 2} & - 2 (n-3) \\
   - {n-1 \choose 2} & - {n-2 \choose 2} & 2n-5 
\end{array} \right )$$
with respect to $a, b$ and $c$.  The left eigenvectors of this are
$$(1,1,1) \hspace{.25in} , \hspace{.25in} 
   \left ( \frac{n-1}{2} \; , \; \frac{n-3}{4} \, , \, -1 \right )
   \hspace{.5in} \mbox{ and } \hspace{.5in} 
   \left ( {n-1 \choose 2} \; , \; - \frac{n-2}{2} \, , \, 1 \right ) \, .$$
The corresponding eigenvalues are 
$$0 \; , \; \frac{2}{3} \frac{(n+1)(n-2)}{n(n-1)} \; \mbox{ and } \; 
   \frac{2}{3} \frac{n-3}{n-1} \; .$$

The equations of mean motion are 
$$\E (W(t+1) - W(t)  \| W(t) - \cc + \ee \hh) 
   = \frac{6}{n} x^{-1} (M - I) \hh + O(|\hh|^2) \, , $$
whence the point $\cc$ is attracting for sufficiently small $x$
if and only if the real parts of all eigenvalues of $M$ are less than~1.
We have identified that this is so, and therefore $\cc$ is attracting.
The hypotheses of Proposition~\ref{pr:multi} are therefore satisfied with
a quadratic Lyapunov function, and Theorem~\ref{th:transience time}
follows from Theorem~\ref{th:rate}.    $\Cox$


\begin{thebibliography}{YMN}

\bibitem[AI97]{AI97}
Anderlini, L. and Ianni, A. (1997).  Learning on a torus. In: {\em
The dynamics of norms}, ed. C. Bicchieri, R. Jeffrey, and B.
Skyrms (Cambridge University Press, Cambridge) 87--107.

\bibitem[BI84]{BI84}
Bannai, E. and Ito, T. (1984).  {\em Algebraic combinatorics I: 
association schemes.}  Benjamin/Cummings: Menlo Park, CA.

\bibitem[BCN89]{BCN89}
Brouwer, A., Cohen, A. and Neumaier, A. (1989).  {\em Distance-regular 
graphs.}  Springer-Verlag: Berlin.

\bibitem[BW00]{BW00}
Barrat, A. and Weigt, M. (2000).  On the properties of small-world
network models.  {\em Europ. Phys. J. B} {\bf 13}, 547.

\bibitem[Ben96]{Ben96}
Bena\"im, M. (1996).  A dynamical system approach to stochastic
approximations.  {\em SIAM J. Control Opt.} {\bf 34} 437 - 472.

\bibitem[BH95]{BH95}
Bena\"im, M. and Hirsch, M. (1995).  Dynamics of Morse-Smale urn
processes.  {\em Ergodic Theory and Dynamical Systems} {\bf 15}
1005 - 1030.

\bibitem[BP03]{BP03}
Benjamini, I. and Pemantle, R. (2003).  Probabilities for cooled 
Brownian motion to linger near the top of a hill, and application 
to a market share model.  {\em Preprint.}

\bibitem[BL02]{BL02}
Bonacich, P. and Liggett, T. (2002).  Asymptotics of a matrix-valued 
Markov chain arising in sociology.  {\em Preprint.}

\bibitem[BM55]{BM55}
Bush, R. and Mosteller, F. (1955).  {\em Stochastic models for learning.}
John Wiley \& Sons: New York.

\bibitem[Dav90]{Dav90}
Davis, B. (1990).  Reinforced random walk.  {\em Prob. Th. Rel. Fields}
{\bf 84}, 203--229.

\bibitem[DZ93]{DZ93}
Dembo, A. and Zeitouni, O. (1993).  {\em Large deviations techniques and 
applications.}  Jones and Bartlett: Boston.

\bibitem[EP23]{EP23}
Eggenberger, F. and P\'olya, G. (1923).  Uber die Statistik verketter
vorg\"age.  {\em Zeit. Angew. Math. Mech.} {\bf 1},~279-289.

\bibitem[Ell93]{Ell93}
Ellison, G. (1993).  Learning, local interaction, and
coordination.  {\em Econometrica} {\bf 61},~1047--1071.

\bibitem[Fr65]{Fr65}
Freedman, D. (1965).  Bernard Friedman's urn.  {\em Ann. Math.
Stat.} {\bf 36} 956 - 970.

\bibitem[Fr49]{Fr49}
Friedman, B. (1949).  A simple urn model.  {\em Comm. Pure Appl. Math.}
{\bf 2},~59--70.

\bibitem[FK93]{FK93}
Fudenberg, D. and Kreps, K. (1993).  Learning mixed equilibria.
{\em Games and Econ. Beh.} {\bf 5},~320--367.

\bibitem[IT69]{IT69}
Iosifescu, M. and Theodorescu, R. (1969).  {\em Random processes
and learning.} Springer-Verlag: New York.

\bibitem[Lak81]{Lak81}
Lakshmivarahan, S. (1981).  {\em Learning algorithms: theory
and applications.} Springer-Verlag: New York.

\bibitem[Lim01]{Lim01}
Limic, V. (2001).  Attracting edge property for a class of reinforced
random walks.  {\em Preprint.} http://www.math.cornell.edu/~limic/

\bibitem[FM96]{FM96}
Flache, A. and Macy, M. (1996).  The weakness of strong ties: collective
action failure in a highly cohesive group.  {\em J. Math. Sociol.}
{\bf 21},~3--28.

\bibitem[MS82]{MS82}
Maynard Smith, J. (1982).  {\em Evolution and the theory of games.}
Cambridge University Press: Cambridge.

\bibitem[Nor72]{Nor72}
Norman, M. (1972).  {\em Markov processes and learning models.}
Academic Press: New York.

\bibitem[Pem90]{Pem90}
Pemantle, R. (1990).  Nonconvergence to unstable points in urn
models and stochastic approximations.  {\em Ann. Probab.} {\bf 18}
698 - 712.

\bibitem[PS03a]{PS03a}
Pemantle, R. and Skyrms, B. (2003).  Network formation by reinforcement 
learning: the long and medium run.  {\em Preprint.}

\bibitem[PS03b]{PS03b}
Pemantle, R. and Skyrms, B. (2003).  Reinforcement schemes may take 
a long time to exhibit limiting behavior.  {\em In Preparation.}

\bibitem[RM51]{RM51}
Robbins, H. and Monro, S. (1951).  A stochastic approximation method.
{\em Ann. Math. Stat.} {\bf 22}, 400--407.

\bibitem[RE95]{RE95}
Roth, A. and Erev, I. (1995).  Learning in extensive form games:
experimental data and simple dynamic models in the intermediate term.
{\em Games and Economic Behavior} {\bf 8},~164--212.


\bibitem[Ter96]{Ter96}
Terwilliger, P. (1996).  Algebraic Combinatorics.  {\em Lecture Notes,
posted January 24, 1996.}

\bibitem[WS98]{WS98}
Watts, D. and Strogatz, S. (1998) Collective dynamics of
``small-world'' networks. {\em Nature} {\bf 393},~440--442.

\end{thebibliography}
\end{document}